\documentclass[11pt]{amsart}
\usepackage{graphicx}

\usepackage{latexsym,amsmath}
\usepackage{amsmath,amsthm}
\usepackage{amsfonts}
\usepackage[psamsfonts]{amssymb}

\theoremstyle{plain} 
\newtheorem{theorem}{Theorem}[section]
\newtheorem{corollary}[theorem]{Corollary}

\newtheorem{lemma}{Lemma}[subsection]
\newtheorem{proposition}[theorem]{Proposition}

\theoremstyle{definition} 
\newtheorem{definition}[theorem]{Definition}
\theoremstyle{definition} 

\theoremstyle{remark} 

\theoremstyle{remark} 
\newtheorem{remark}[theorem]{Remark}
\newtheorem*{remark*}{Remark}
\numberwithin{equation}{section}
 
\newcommand{\beqa}{\begin{eqnarray}}
\newcommand{\eeqa}{\end{eqnarray}}
 
\newcommand{\bseq}{\begin{subequations}}
 
\newcommand{\eseq}{\end{subequations}}

\newcommand{\lp}{\left(}
\newcommand{\rp}{\right)}

\newcommand{\al}{\alpha}
\newcommand{\g}{\gamma}
\newcommand{\si}{\sigma}
\newcommand{\la}{\lambda}

\renewcommand{\Psi}{\overline{\Phi}}

\newcommand{\uu}{\nearrow}
\newcommand{\se}{\searrow}

\newcommand{\opt}{\textsf{opt}}

\renewcommand{\P}{\mathsf{P}} 

\newcommand{\E}{\mathsf{E}}
\newcommand{\Var}{\mathsf{Var}}

\newcommand{\R}{\mathbb{R}}

\newcommand{\F}{\mathcal{F}}

\newcommand{\G}{\mathcal{F}}

\renewcommand{\d}{\mathrm{d}}

\newcommand{\vp}{\varepsilon}

\begin{document}

\begin{abstract}
Let $(S_0,S_1,\dots)$ be a supermartingale relative to a nondecreasing sequence of $\sigma$-algebras $(H_{\le0},H_{\le1},\dots)$, with $S_0\le0$ almost surely (a.s.) and differences $X_i:=S_i-S_{i-1}$. Suppose that for every $i=1,2,\dots$ there exist $H_{\le(i-1)}$-measurable r.v.'s $C_{i-1}$ and $D_{i-1}$ and a positive real number $s_i$ such that 
$C_{i-1}\le X_i\le D_{i-1}$ and 
$D_{i-1}-C_{i-1}\le 2 s_i$
a.s. Then for all real $t$ and natural $n$
\begin{equation*}
\E f_t(S_n)\le\E f_t(sZ),
\end{equation*}
where $f_t(x):=\max(0,x-t)^5$,
$s:=\sqrt{s_1^2+\dots+s_n^2}$,
and $Z\sim N(0,1)$. In particular, this implies
\begin{equation*}
\P(S_n\ge x)\le c_{5,0}\P(Z\ge x/s)\quad\forall x\in\R,
\end{equation*}
where $c_{5,0}=5!(e/5)^5=5.699\dots.$ Results for $\max_{0\le k\le n}S_k$ in place of $S_n$ and for concentration of measure also follow.
\end{abstract}

\title[Normal domination of (super)martingales]{{\Large On normal domination of (super)martingales}}
\date{\today}
\author{Iosif Pinelis}
\address{ Department of Mathematical Sciences\\
Michigan Technological University\\
Houghton, Michigan 49931 }
\email{ipinelis@math.mtu.edu}
\maketitle

\setcounter{tocdepth}{4}
\tableofcontents


\makeatletter
\makeatother

\section{Introduction}\label{intro}

The sharp form, 
\begin{equation}\label{eq:khin}
\E f\lp\vp_1a_1+\dots+\vp_n a_n\rp\le\E f(Z),
\end{equation}
of Khinchin's inequality for $f(x)=|x|^p$ 
for the normalized Rademacher sum $\vp_1a_1+\dots+\vp_n a_n$, with
$$a_1^2+\dots+a_n^2=1,$$
was proved by
Whittle (1960) \cite{whittle} for $p\ge3$ and Haagerup (1982) \cite{haag} for $p\ge2$; here and elsewhere, 
the $\vp_i$'s are independent Rademacher random variables (r.v.'s), so that $\P(\vp_i=1)=\P(\vp_i=-1)=1/2$ for all $i$, and $Z\sim N(0,1)$. 

For $f(x)=e^{\la x}$ ($\la\ge0$), this inequality follows from Hoeffding (1963) \cite{hoeff}, whence
$$
\P\lp\vp_1a_1+\dots+\vp_n a_n\ge x\rp\le 
\inf_{\la\ge0}\frac{\E e^{\la Z}}{e^{\la x}}
=e^{-x^2/2},\quad x\ge0.
$$
Since 
$\P(Z\ge x)\sim\frac1{x\sqrt{2\pi}}e^{-x^2/2}$ $(x\to\infty)$,
a factor $\asymp\frac1x$ is ``missing" here. 
The apparent cause of this deficiency is that the class of the exponential moment functions $f(x)=e^{\la x}$ ($\la\ge0$) is too small (and so is the class of the power functions $f(x)=|x|^p$). 
 
Consider the much richer classes of functions
$\G_+^{(\al)}$ ($\al\ge0$), consisting of all the functions $f\colon\R\to\R$ given by the formula
$$
f(x)=\int_{-\infty}^\infty (x-t)_+^\al\,\mu(dt),\ u\in\R,
$$
where $\mu\ge0$ is a Borel measure,
$x_+:=\max(0,x)$, $x_+^\al:=(x_+)^\al$,
$0^0:=0$.

It is easy to see \cite[Proposition 1(ii)]{pin99} that
\begin{equation}\label{eq:F-al-beta}
0\le\beta<\al\quad\text{implies}\quad\G_+^{(\al)}\subseteq\G_+^{(\beta)}.
\end{equation}

\begin{proposition}\label{prop:F-al}
\cite{binom}\ 
For natural $\al$, one has $f\in\G_+^{(\al)}$ if and only if $f$ has finite derivatives $f^{(0)}:=f,f^{(1)}:=f',\dots,f^{(\al-1)}$ on $\R$ such that $f^{(j)}(-\infty)=0$ for $j=0,1,\dots,\al-1$ and $f^{(\al-1)}$ is convex on $\R$. 
\end{proposition}

It follows from Proposition \ref{prop:F-al} that, for every $t\in\R$, every $\beta\ge\al$, and every $\la>0$, the functions $u\mapsto(u-t)_+^\beta$ and $u\mapsto e^{\la(u-t)}$
belong to $\G_+^{(\al)}$, while
the functions $u\mapsto|u-t|^\beta$ and $u\mapsto\cosh\la(u-t)$
belong to $\G^{(\al)}$.  

Eaton (1970) \cite{eaton1} proved the Khinchin-Whittle-Haagerup inequality \eqref{eq:khin} for a class of moment functions, which essentially coincides with the class $\G_+^{(3)}$.
Based on asymptotics, numerics, and a certain related inequality, 
Eaton (1974) \cite{eaton2} conjectured that the mentioned moment comparison inequality of his implies that 
$$\P\,\lp\vp_1a_1+\dots+\vp_n a_n\ge x\rp\le \frac{2e^3}9\,\frac1{x\sqrt{2\pi}}e^{-x^2/2}\quad \forall x>\sqrt2.
$$
Pinelis (1994) \cite{pin94} proved the following improvement of this conjecture: 
\begin{equation}\label{eq:pin94}
\P\,\lp\vp_1a_1+\dots+\vp_n a_n\ge x\rp\le\frac{2e^3}9\,\P(Z\ge x)\quad
\forall x\in\R, 
\end{equation}
as well as certain multidimensional extensions of these results.

Later it was realized in Pinelis (1998) \cite{pin98} that the reason why it is possible to extract tail comparison inequality \eqref{eq:pin94} from the Khinchin-Eaton moment comparison inequality \eqref{eq:khin} for $f\in\G_+^{(3)}$ is that the tail function $x\mapsto\P(Z\ge x)$ is log-concave.
This realization resulted in a general device, which allows one to extract the optimal tail comparison inequality from an appropriate moment comparison inequality. The following is a special case of Theorem 4 of Pinelis (1999) \cite{pin99}; see also Theorem 3.11 of Pinelis (1998) \cite{pin98}.
\begin{theorem}
\label{th:comparison} 
Suppose that $0\le\beta\le\al$, $\xi$ and $\eta$ are real-valued r.v.'s, and the tail function $u\mapsto\P(\eta\ge u)$ is log-concave on $\R$. Then the comparison inequality 
\begin{equation}\label{eq:comp-al}
\E f(\xi)\le\E f(\eta)\quad\text{for all }f\in\G_+^{(\al)}
\end{equation}
implies
\begin{equation}\label{eq:comp-beta}
\E f(\xi)\le c_{\al,\beta}\,\E f(\eta)\quad\text{for all }f\in\G_+^{(\beta)}
\end{equation}
and, in particular, for all real $x$,
\begin{align}
\P(\xi\ge x) &\le \inf_{f\in\G_+^{(\al)}}\,\frac{\E f(\eta)}{f(x)} \label{eq:comp-prob1} \\
&= B_{\opt}(x):=\inf_{t\in(-\infty,x)}\,\frac{\E(\eta-t)_+^\al}{(x-t)^\al} 
\label{eq:comp-prob2} \\
&\le \min\lp c_{\al,0}\,\P(\eta\ge x),\; \inf_{h>0}\,e^{-hx}\,\E e^{h\eta} \rp,
\label{eq:comp-prob3} 
\end{align}
where 
\begin{equation}\label{eq:c(al,beta)}
c_{\al,\beta}:=\frac{\Gamma(\al+1)(e/\al)^\al}{\Gamma(\beta+1)(e/\beta)^\beta}.
\end{equation}
Moreover, the constant $c_{\al,\beta}$ is the best possible in \eqref{eq:comp-beta} and \eqref{eq:c(al,beta)}. 
\end{theorem}

A similar result for the case when $\al=1$ and $\beta=0$ 
is contained in the book by Shorack and Wellner (1986) \cite{shorack-wellner}, pages 797--799. 

\begin{remark}
\label{comparison-remark}
As folows from \cite[Remark 3.13]{pin98}, a useful point is that the requirement of the log-concavity of the tail function $q(u):=\P(\eta\ge u)$ in Theorem \ref{th:comparison}
can be relaxed by replacing $q(x)=\P(\eta\ge x)$ by any [e.g., the least] log-concave majorant of $q$. 
However, then the optimality of $c(\al,\beta)$ is then not guaranteed. 
\end{remark}

Note that 
$c_{3,0}=2e^3/9$, 
which is the constant factor in \eqref{eq:pin94}.
Bobkov, G\"{o}tze, and Houdr\'{e} (2001) \cite{houdre} obtained a simpler proof of inequality \eqref{eq:pin94},
but with a constant factor $12.0099\ldots$ in place of $2e^3/9=4.4634\ldots$. 

Pinelis (1999) \cite{pin99} obtained the ``discrete" improvement of \eqref{eq:pin94}: 
\begin{equation}\label{eq:pin99}
\P\lp\vp_1a_1+\dots+\vp_n a_n\ge x\rp 
\le\frac{2e^3}9\,
\P\lp\frac1{\sqrt n}(\vp_1+\dots+\vp_n)\ge x\rp 
\end{equation}
for all values $x$ of r.v. $\frac1{\sqrt n}(\vp_1+\dots+\vp_n)$.



\section{Domination by normal moments and tails}\label{normal}

\begin{theorem}\label{th:bernoulli} 
Let $S_0\le0,S_1,\dots$ be a supermartingale, with increments $X_i:=S_i-S_{i-1}$, $i=1,2,\dots$. Suppose that for every $i=1,2,\dots$ there exist $H_{\le(i-1)}$-measurable r.v.'s $C_{i-1}$ and $D_{i-1}$ and a positive real number $s_i$ such that 
\begin{gather}
C_{i-1}\le X_i\le D_{i-1} \quad\text{ and } \label{eq:bern-cond1} \\
D_{i-1}-C_{i-1}\le 2 s_i \label{eq:bern-cond2}
\end{gather}
with probability 1. Then for all $f\in\G_+^{(5)}$ and all $n=1,2,\dots$
\begin{equation}\label{eq:bernoulli}
\E f(S_n)\le\E f(sZ),
\end{equation}
where 
$$s:=\sqrt{s_1^2+\dots+s_n^2}$$
and $Z\sim N(0,1)$. 
\end{theorem}

The proof of this and other statements (whenever necessary) are deferred to Section \ref{proofs}. 

By virtue of Theorem \ref{th:comparison}, one has the following corollary under the conditions of Theorem \ref{th:bernoulli}.

\begin{corollary}\label{cor:bern-beta}
For all $\beta\in[0,5]$, all $f\in\G_+^{(\beta)}$, and all $n=0,1,\dots$
\begin{equation}\label{eq:bern-beta}
\E f(S_n)\le c_{5,\beta}\,\E f(sZ).
\end{equation}
In particular, for all real $x$,
\begin{align}
\P(S_n\ge x) &\le \inf_{f\in\G_+^{(5)}}\,\frac{\E f(sZ)}{f(x)} \label{eq:bern-prob1} \\
&= \inf_{t\in(-\infty,x)}\,\frac{\E(sZ-t)_+^\al}{(x-t)^\al} 
\label{eq:bern-prob2} \\
&\le \min\lp c_{5,0}\,\P(sZ\ge x), \inf_{h>0}\,e^{-hx}\,\E e^{hsZ} \rp 
\label{eq:bern-prob3} \\
&= \min\lp c_{5,0}\,\Psi\lp\frac{x}{s}\rp, \exp\lp -\frac{x^2}{2s^2} \rp\rp, 
\label{eq:bern-prob4} 
\end{align}
and 
$$c_{5,0}=5!(e/5)^5=5.699\dots.$$
\end{corollary}

The upper bound $\exp\left(-\frac{x^2}{2s^2}\right)$ was obtained by Hoeffding (1963) \cite{hoeff} for the case when the $C_{i-1}$'s and $D_{i-1}$'s are non-random.

The upper bound \eqref{eq:bern-prob3} -- but with constant factor $435$ in place of $c_{5,0}=5.699\dots$ -- was obtained in \cite{bent-isr} for the case when $(S_i)$ is a martingale. 


\begin{theorem}
\label{th:bernoulli-improved} Let $S_0\le0,S_1,\dots$ be a supermartingale, with increments $X_i:=S_i-S_{i-1}$, $i=1,2,\dots$. Suppose that for every $i=1,2,\dots$ there exist a positive $H_{\le(i-1)}$-measurable r.v. $D_{i-1}$ and a positive real number $s_i$ such that 
\begin{gather}
X_i\le D_{i-1} \quad \text{and} \label{eq:bern-improv-cond1} \\ \frac12\left(D_{i-1}+\frac{\Var_{i-1}X_i}{D_{i-1}}\right)\le \hat s_i \label{eq:bern-improv-cond2} 
\end{gather}
with probability 1. 
Let 
\begin{equation}\label{eq:hat-s}
\hat s:=\sqrt{\hat s_1^2+\dots+\hat s_n^2}.
\end{equation}
Then one has all the inequalities \eqref{eq:bernoulli}--\eqref{eq:bern-prob4}, only with $s$ replaced by $\hat s$. 
\end{theorem}

\begin{remark}
Theorem~\ref{th:bernoulli} may be considered as a special case Theorem~\ref{th:bernoulli-improved}. Indeed, it can be seen from the proofs of these two theorems (see Lemma~\ref{lem:reduction1} and 
Lemma~3.1 in \cite{binom}),
one may assume without loss of generality that the supermartingales $(S_i)$ in Theorem \ref{th:bernoulli} and \ref{th:bernoulli-improved} are actually martingales with $S_0=0$. Therefore, to deduce Theorem \ref{th:bernoulli} from Theorem \ref{th:bernoulli-improved}, it is enough to observe that for any r.v. $X$ and constants $c<0$ and $d>0$, one has the following implication:
\begin{equation}\label{eq:implic}
\E X=0\ \&\ \P(c\le X\le d)=1\implies\Var X\le|c|d.
\end{equation}
In turn, implication \eqref{eq:implic} follows from \cite{karr}, which reduces the sitation to that of a r.v. $X$ taking on onlyt two values. Alternatively, in light of the duality result \cite[(4)]{pin98}, it is easy to give a direct proof of \eqref{eq:implic}. Indeed, $\E X=0$ and $\P(c\le X\le d)=1$ imply
$$0\ge\E(X-c)(X-d)=\E X^2+cd=\Var X-|c|d.$$

However, rather than deducing Theorem \ref{th:bernoulli} from Theorem \ref{th:bernoulli-improved}, we shall go in the opposite direction, proving Theorem \ref{th:bernoulli-improved} based on Theorem \ref{th:bernoulli}. 

Thus, Theorem \ref{th:bernoulli} is seen as the main result of this paper.
\end{remark}

\begin{remark}\label{rem:bernoulli-improved}
The set of conditions \eqref{eq:bern-improv-cond1}--\eqref{eq:bern-improv-cond2} is equivalent to
$$X_i\le D_{i-1} \quad \text{and}\quad \si_*(D_{i-1}, \E_{i-1}X_i^2) \le s_i
$$
with probability 1, where
\begin{multline*}
\si_*(d_0,\si^2):=
\frac12\,\inf_{d\ge d_0} \left(d+\frac{ \si^2 }d \right)
=\min\left( \si\vee d_0, \frac12\left(d_0+\frac{ \si^2 }{d_0}\right)\right) \\
=\begin{cases}
\si & \text{if } \si\ge d_0, \\
\frac12\, \left(d+\frac{ \si^2 }d \right) & \text{if } \si<d_0, 
\end{cases}
\end{multline*}
for positive $\si$ and $d_0$. 
This follows simply because the inequalities $X_i\le D_{i-1}$ and $d\ge D_{i-1}$ imply $X_i\le d$.  
\end{remark}

From the ``right-tail" bounds stated above, ``two-tail" ones immediately follow:

\begin{corollary}
\label{cor:bernoulli1} Let $S_0=0,S_1,\dots$ be a martingale, with increments $X_i:=S_i-S_{i-1}$, $i=1,2,\dots$. Suppose that conditions \eqref{eq:bern-cond1} and \eqref{eq:bern-cond2} hold. Then inequalities \eqref{eq:bernoulli} and \eqref{eq:bern-beta} hold for all $f\in\G^{(5)}$ and $f\in\G^{(\beta)}$ ($\beta\in[0,5]$), rather than only for all $f\in\G_+^{(5)}$ and $f\in\G_+^{(\beta)}$, respectively.  
\end{corollary}

\begin{corollary}
\label{cor:bernoulli2} Let $S_0=0,S_1,\dots$ be a martingale, with increments $X_i:=S_i-S_{i-1}$, $i=1,2,\dots$. Suppose that condition \eqref{eq:bern-improv-cond2} holds, and condition \eqref{eq:bern-improv-cond1} holds for $|X_i|$ in place of $X_i$. Then inequalities \eqref{eq:bernoulli} and \eqref{eq:bern-beta} with $s$ replaced by $\hat s$ hold for all $f\in\G^{(5)}$ and $f\in\G^{(\beta)}$ ($\beta\in[0,5]$), rather than only for all $f\in\G_+^{(5)}$ and $f\in\G_+^{(\beta)}$, respectively.  
\end{corollary}

That $(S_0,S_1,\dots)$ in Theorems \ref{th:bernoulli} and \ref{th:bernoulli-improved} is allowed to be a supermartingale (rather than only a martingale) makes it convenient to use the simple but powerful truncation tool. (Such a tool was used, for example, in \cite{pin85} to prove limit theorems for large deviation probabilities based only on precise enough probability inequalities and without using Cram\'{e}r's transform, the standard device in the theory of large deviations.) Thus, for instance, one has the following corollary from Theorem \ref{th:bernoulli-improved}.

\begin{corollary}\label{cor-truncate} 
Let $S_0\le0,S_1,\dots$ be a supermartingale, with increments $X_i:=S_i-S_{i-1}$, $i=1,2,\dots$. For every $i=1,2,\dots$, let $D_{i-1}$ be a positive $H_{\le(i-1)}$-measurable r.v. and let $s_i$ be a positive real number
such that \eqref{eq:bern-improv-cond2} holds (while \eqref{eq:bern-improv-cond1} does not have to). Let $\hat s$ be still defined by \eqref{eq:hat-s}. 

Then for all real $x$
\begin{align}
\P(S_n\ge x) &\le \P\lp\max_{1\le i\le n}\,\frac{X_i}{D_{i-1}} \ge 1 \rp +
\min\lp c_{5,0}\,\Psi\lp\frac{x}{s}\rp, \exp\lp -\frac{x^2}{2s^2} \rp\rp 
\label{eq:truncate1} \\ 
			&\le \sum_{1\le i\le n}\P\lp X_i\ge D_{i-1} \rp +
\min\lp c_{5,0}\,\Psi\lp\frac{x}{s}\rp, \exp\lp -\frac{x^2}{2s^2} \rp\rp .
\label{eq:truncate2}  
\end{align}
\end{corollary}

These bounds are much more precise than the exponential bounds in \cite{fuk-nagaev,fuk,nag}. 



\section{Maximal inequalities}\label{max}

Introduce 
$$M_n:=\max_{0\le k\le n}S_k.$$

\begin{theorem}\label{th:max}
Let $(S_0=0,S_1,\dots)$ be a martingale. Then the upper bounds on $\P(S_n\ge x)$ given in Corollary \ref{cor:bern-beta} and Theorem \ref{th:bernoulli-improved} are also upper bounds on $\P(M_n\ge x)$, under the same conditions: \eqref{eq:bern-cond1}-\eqref{eq:bern-cond2} and \eqref{eq:bern-improv-cond1}-\eqref{eq:bern-improv-cond2}, respectively.  
\end{theorem}

\begin{theorem}\label{th:al-beta-Mn}
Let $0\le\beta\le\al$ and $x>t$, and let $(S_n)$ be a martingale or, more generally, a submartingale. 
Assume, moreover, that
$\al>1$. 
Then, for any natural $n$,
\begin{equation}\label{eq:al-beta-Mn}
\E(M_n-x)_+^\beta\le k_{1;\al,\beta}\frac{\E(S_n-t)_+^\al}{(x-t)^{\al-\beta}},
\end{equation}
where
\begin{equation}\label{eq:k1}
k_{1;\al,\beta}:=
\sup_{\si>0} \si^{-\beta(\al-1)}
\left(\int_0^\si \frac{\beta s^{\beta-1}\,ds}{1+s}\right)^\al
\end{equation}
if $\beta>0$, and $k_1(\al,0):=1$. 
The particular cases of \eqref{eq:al-beta-Mn}, corresponding to $\beta=0$ and $\beta=\al$, respectively, are Doob's 
inequalities
\begin{equation}\label{eq:doob-0}
\P(M_n\ge x)\le\frac{\E(S_n-t)_+^\al}{(x-t)^\al}
\end{equation}
and
\begin{equation}\label{eq:doob-al}
\E(M_n)_+^\al\le\lp\frac\al{\al-1}\rp^\al\,\E(S_n)_+^\al.
\end{equation}
\end{theorem}

\begin{theorem}\label{th:max-beta}
Let $(S_0=0,S_1,\dots)$ be a martingale. Then inequalities \eqref{eq:bern-beta} and \eqref{eq:bern-prob1} hold if $S_n$ is replaced there by $M_n$ and $c_{5,\beta}$ by $\dfrac{k_{1;\al,\beta}}{k_{\al,\beta}}\,c_{5,\beta}$, under the same conditions: \eqref{eq:bern-cond1}-\eqref{eq:bern-cond2} and \eqref{eq:bern-improv-cond1}-\eqref{eq:bern-improv-cond1}, respectively.  
\end{theorem}

Similarly, results of \cite{binom} can be extended.

\begin{remark}
Note that 
$$\int_0^\si \frac{\beta s^{\beta-1}\,ds}{1+s}=\si^\beta
{}_2\!F_1(\beta,1;1+\beta;-\si)
=\beta\int_0^1(1-u)^{\beta-1}(1+\si u)^{-\beta}\,du,$$
where ${}_2\!F_1$ is a hypergeometric function.  
Note also that there is some $\si_{\al,\beta}\in(0,\infty)$ such that the expression under the $\sup$ sign in \eqref{eq:k1} is increasing in $\si\in(0,\si_{\al,\beta})$ and decreasing in $\si\in(\si_{\al,\beta},\infty)$; this can be seen from the proof of Proposition \ref{prop:k1-k3}. Thus, the $\sup$ is attained at the unique point $\si_{\al,\beta}$. 
\end{remark}

\begin{proposition}\label{prop:k1-k2}
Let $\al$ and $\beta$ be as in Theorem \ref{th:al-beta-Mn}.
Then
\begin{equation}\label{eq:k1-k2}
k_{1;\al,\beta}\le k_{2;\al,\beta}:=\frac{\Gamma(1+\beta)\Gamma(\al-\beta)}{\Gamma(\al)}.
\end{equation}
\end{proposition}

\begin{remark}
$$k_2(\al,0)=k(\al,0)=1=k_1(\al,0).$$
\end{remark}

\begin{proposition}\label{prop:al-beta}
Let $0\le\beta<\al$, $x>t$, and 
\begin{equation}\label{eq:k}
k_{\al,\beta}:=\frac{\beta^\beta(\al-\beta)^{\al-\beta}}{\al^\al}.
\end{equation}
Then
\begin{equation}\label{eq:al-beta}
\forall u\in\R\quad(u-x)_+^\beta\le k_{\al,\beta}\frac{(u-t)_+^\al}{(x-t)^{\al-\beta}},
\end{equation}
and $k_{\al,\beta}$ is the best constant here. (The values at $\beta=0$ are understood here as the corresponding limits as $\beta\downarrow0$.) 
\end{proposition}

\begin{proposition}\label{prop:al-beta-Sn}
Let $0\le\beta\le\al$ and $x>t$, and let $(S_n)$ be a martingale or, more generally, a submartingale. Then, for any natural $n$,
\begin{equation}\label{eq:al-beta-Sn}
\E(S_n-x)_+^\beta\le k_{\al,\beta}\frac{\E(S_n-t)_+^\al}{(x-t)^{\al-\beta}},
\end{equation}
and $k_{\al,\beta}$ is the best constant here. 
\end{proposition}

\begin{proposition}\label{prop:k1-k3}
Let $\al$ and $\beta$ be as in Theorem \ref{th:al-beta-Mn}.
Then
\begin{equation}\label{eq:k1-k3}
k_{1;\al,\beta}\le k_{3;\al,\beta}:=k_{\al,\beta}\left(\frac\al{\al-1}\right)^\al,
\end{equation}
where $k_{\al,\beta}$ is defined by \eqref{eq:k}.
\end{proposition}

\begin{proposition}\label{prop:k1=k3}
Let $\al>1$.
Then
\begin{equation}\label{eq:k1=k3}
k_1(\al,\al)=k_3(\al,\al)=\left(\frac\al{\al-1}\right)^\al.
\end{equation}
\end{proposition}

\begin{corollary}\label{cor:k,k1,k2,k3}
Let $\al$ and $\beta$ be as in Theorem \ref{th:al-beta-Mn}.
Then
\begin{equation}\label{eq:k,k1,k2,k3}
k(\al,\al)\le k_1(\al,\al)\le k_2(\al,\al)\wedge k_3(\al,\al);
\end{equation}
at that
\begin{equation}\label{eq:k,k1,k2,k3(0)}
k(\al,0)=k_1(\al,0)=k_2(\al,0)=1,
\end{equation}
while
\begin{equation}\label{eq:k,k1,k2,k3(al)}
k_1(\al,\al)=k_3(\al,\al)=\left(\frac\al{\al-1}\right)^\al>k(\al,\al)=1.
\end{equation}
\end{corollary}

\section{Concentration inequalities for separately Lipschitz functions}\label{concentr}

\begin{definition}\label{def:lip}
Let us say that a real-valued function $g$ of $n$ (not necessarily real-valued) arguments is {\em separately Lipschitz} if it satisfies a Lipschitz type condition in each of its arguments:
\begin{equation}\label{eq:Lip}
|g(x_1,\dots,x_{i-1},\tilde x_i,x_{i+1},\dots,x_n) -
g(x_1,\dots,x_n)| \le \rho_i(\tilde x_i,x_i) < \infty
\end{equation}
for all $i$ and all $x_1,\dots,x_n,\tilde x_i$, where $\rho_i(\tilde x_i,x_i)$ depends only on $\tilde x_i$ and $x_i$.
Let the {\em radius} of the separately Lipschitz function $g$ be defined as 
$$r := \sqrt{ r_1^2 +\dots+ r_n^2},$$
where
\begin{equation}\label{eq:r_i}
r_i:=\frac12\,\sup_{\tilde x_i,x_i} \rho_i(\tilde x_i,x_i).
\end{equation}
\end{definition}

The concentration inequalities given in this section follow from martingale inequalities given in Section~\ref{normal}. 
Their proofs here are based on the improvements given in \cite{pin81a} and \cite{pin-sakh} of the method of Yurinskii (1974) \cite{yurin}; cf. \cite{mcdiarmid1,mcdiarmid2} and \cite{bent-isr}. 

Papers \cite{yurin}, \cite{pin81a}, and \cite{pin-sakh} deal mainly with separately Lipschitz function $g$ of the form
$$g(x_1,\dots,x_n)=\|x_1+\dots+x_n\|,$$
where the $x_i$'s are vectors in a normed space; 
however, it was already understood there that the methods would work for much more general functions $g$ -- see \cite[Remark 1]{pin-sakh}. 
In a similar fashion, various concentration inequalities for general functions $g$ were obtained in 
\cite{mcdiarmid1,mcdiarmid2} and 
\cite{bent-isr}.


\begin{theorem}\label{th:concentr}
Suppose that a r.v. $Y$ can be represented as a real-valued function $g$ of independent 
(not necessarily real-valued) r.v.'s $X_1,\dots,X_n$:
$$Y=g(X_1,\dots,X_n),$$
where $g$ is separately Lipschitz with radius $r$. 
Then
\begin{align} 
\E f(Y - \E Y) & \le \E f(rZ)\quad\text{for all }f\in\G^{(5)}\quad\text{and} \label{eq:concentr5} \\
\E f(Y - \E Y) & \le c_{5,\beta}\,\E f(rZ)\quad\text{for all }\beta\in[0,5]\text{ and all }f\in\G^{(\beta)}, 
\label{eq:concentrBeta}
\end{align}
where $Z\sim N(0,1)$. 
In particular, for all real $x$,
\begin{equation} 
\P(Y - \E Y \ge x) \le c_{5,0}\,\P(rZ\ge x) = c_{5,0}\,\Psi\lp\frac{x}{r}\rp.  
\label{eq:concentrTail}
\end{equation}
\end{theorem}

\begin{proposition}\label{prop:LipRelaxed}
Inequalities \eqref{eq:concentr5}, \eqref{eq:concentrBeta}, and \eqref{eq:concentrTail} will hold if the conditions of
Theorem \ref{th:concentr}
are relaxed so that $r_i$ is replaced by
\begin{multline}\label{eq:LipRelaxed}
\hat r_i:=\frac12\,\sup_{x_1,\dots,x_i,\tilde x_i}|\E g(x_1,\dots,x_{i-1},\tilde x_i,X_{i+1},\dots,X_n) 
\\ - 
\E g(x_1,\dots,x_i,X_{i+1},\dots,X_n)|,
\end{multline}
for every $i$.
Note that $\hat r_i\le r_i$ for all $i$. 
\end{proposition}

\begin{remark}\label{rem:tighter}
The upper bound given by \eqref{eq:concentrTail} can be replaced by the tighter bound
$$\min\left( \exp\left(-\frac{x^2}{2r^2}\right), c_{5,0}\,\Psi\lp\frac{x}{r}\rp\right),$$
which is less than $\exp\left(-\frac{x^2}{2r^2}\right)$ for all $\frac{x}{r}\ge1.89$.
\end{remark}

The foregoing conditions can be modified as follows. 

\begin{theorem}\label{th:concentr-modif}
Suppose that
\begin{align}
\Xi_i(x_1,\dots,x_{i-1},x_i)&:=\E g(x_1,\dots,x_{i-1},x_i,X_{i+1},\dots,X_n) \label{eq:Xi_i}
\\ & - 
\E g(x_1,\dots,x_{i-1},X_i,X_{i+1},\dots,X_n) \notag \\
& \le D_{i-1}(x_1,\dots,x_{i-1}), \label{eq:LipModif}
\end{align}
and 
\begin{equation}\label{eq:Lip-s_i}
\frac12\, \left( D_{i-1}(x_1,\dots,x_{i-1}) + 
\frac{ \E\Xi_i(x_1,\dots,x_{i-1},X_i)^2 }{ D_{i-1}(x_1,\dots,x_{i-1}) }
\right)
\le s_i,
\end{equation}
for all $i$ and all $x_1,\dots,x_{i-1},x_i$, where $D_{i-1}>0$ depends only on $i$ and $x_1,\dots,x_{i-1}$, and $s_i$ depends only on $i$.
Let 
$$s:=\sqrt{s_1^2+\dots+s_n^2}.$$
Then inequalities \eqref{eq:concentr5}, \eqref{eq:concentrBeta}, and \eqref{eq:concentrTail} will hold if $r$ is replaced there by $s$. 
\end{theorem}

The next two propositions show how to obtain good upper bounds on $\Xi_i(x_1,\dots,x_{i-1},x_i)$ and $\E\Xi_i(x_1,\dots,x_{i-1},X_i)^2$, to be used in Theorem \ref{th:concentr-modif}.

\begin{proposition}\label{prop:sigma-xi-bound}
If $g$ is separately Lipschitz so that \eqref{eq:Lip} holds, then for all $i$ and all $x_1,\dots,x_{i-1}$, 
\begin{equation}\label{eq:sigma-xi-bound}
\E\Xi_i(x_1,\dots,x_{i-1},X_i)^2 \le \inf_{x_i}\E\rho_i(X_i,x_i)^2 \le \E\rho_i(X_i,\E X_i)^2.
\end{equation}
If, moreover, the function $g$ is convex in each of its arguments, then for all $i$ and all $x_1,\dots,x_i$, 
\begin{equation}\label{eq:sigma-bound}
\Xi_i(x_1,\dots,x_{i-1},x_i) \le 
\rho_i(x_i,\E X_i).
\end{equation}
\end{proposition}

\begin{remark}
We do not require that $\rho_i$ be a metric. However, the smallest possible $\rho_i$, which is the supremum of the left-hand side of \eqref{eq:Lip} over all $x_1,\dots,x_{i-1},x_{i+1},\dots,x_n$, is necessarily a metric. 
Note also that, for $r_i$ defined by \eqref{eq:Lip},
$$\rho_i(x_i,\E X_i)=\rho_i(x_i,0)\le\frac12\,r_i$$
for all $x_i$, provided the following conditions: (i) $\rho_i$ is the smallest possible and, moreover, 
is a norm; (ii) $X_i$ is symmetrically distributed; and (iii) $x_i$ belongs to the support of the distribution of $X_i$.
\end{remark}

\begin{corollary}\label{cor:banach}
Let here $X_1,\dots,X_n$ be independent r.v.'s with values in a separable Banach space with norm $\|\cdot\|$, and let 
$$Y:= \| X_1+\dots+X_n \|.$$
Suppose that, with probability 1,
\begin{align}
\|X_i-\E X_i\| \le d_i \label{eq:banach-d_i}
\end{align}
and 
\begin{equation}\label{eq:banach-s_i}
\frac12\, \left( d_i + 
\frac{ \E\|X_i-\E X_i\|^2 }{ d_i }
\right)
\le s_i,
\end{equation}
for all $i$, where $d_i>0$ and $s_i>0$ are non-random constants.
Let 
$$s:=\sqrt{s_1^2+\dots+s_n^2}.$$
Then inequalities \eqref{eq:concentr5}, \eqref{eq:concentrBeta}, and \eqref{eq:concentrTail} will hold if $r$ is replaced there by $s$. 
\end{corollary}


\section{Proofs}\label{proofs}

\subsection{Proofs for Section \ref{normal}} \label{subsec:normal}

Let us first observe that Theorem~\ref{th:bernoulli} 
can be easily reduced to the case when $(S_n)$ is a martingale. This is implied by the following two lemmas. 

The next lemma is obvious and stated here for the convenience of reference.

\begin{lemma}\label{lem:reduction1}
Let $(S_n)$ be a supermartingale as in Theorem \ref{th:bernoulli}, 
so that 
conditions \eqref{eq:bern-cond1} and \eqref{eq:bern-cond2} 
are satisfied. 
Let 
\begin{multline*}
\tilde X_i:=X_i-\E_{i-1}X_i,\quad
\tilde C_{i-1}:=C_{i-1}-\E_{i-1}X_i,\quad\text{and}\quad
\tilde D_{i-1}:=D_{i-1}-\E_{i-1}X_i.
\end{multline*}
Then $\tilde X_i$ is $H_{\le i}$-measurable,
$\tilde C_{i-1}$ and $\tilde D_{i-1}$ are $H_{\le(i-1)}$-measurable, and one has
\begin{gather*}
X_i \le \tilde X_i,\\
\E_{i-1}\tilde X_i=0,\\
\tilde C_{i-1}\le \tilde X_i\le \tilde D_{i-1}, \quad\text{ and }  \\
\tilde D_{i-1}-\tilde C_{i-1}\le 2 s_i 
\end{gather*}
with probability 1. 
\end{lemma}

\begin{proof}[Proof of Theorem \ref{th:bernoulli}] is similar to the proof of 
Theorem~2.1 in\cite{binom} but based on the following lemma, in place of Lemma~3.2 in\cite{binom}.
(Also, one has to refer here to Lemma \ref{lem:reduction1} instead of Lemma~3.1 in\cite{binom}.)
\end{proof}

\begin{lemma}
\label{lem:bernoulli} Let $X$ be a r.v. such that $\mathsf{E} X=0$ and $%
c\le X\le d$ with probability 1 for some real constants $c$ and $d$ (whence $%
c\le0$ and $d\ge0$). Let $Z\sim N(0,1)$. Then for all $f\in\G_+^{(5)}$
\begin{equation}  \label{eq:bern-moment-compar}
\E f(X)\le\E f((d-c)Z).
\end{equation}
\end{lemma}

\begin{proof}
This proof is rather long.
Let $\mathcal{X}_{c,d}$ be the set of all r.v.'s $X$ such that $\mathsf{E}X=0
$ and $c\leq X\leq d$ with probability 1. In view of \cite{karr} (say), for any given real $t$, a maximum of $\mathsf{E}%
f_{t}(X)$ over all r.v.'s $X$ in $\mathcal{X}_{c,d}$ is attained when $X$
takes on only two values, say $a$ and $b$, in the interval $[c,d]$. Since
the function $f_{t}$ is convex, it then follows that, without loss of
generality (w.l.o.g.), $a=c$ and $b=d$. $\Bigl($ Indeed, $\mathsf{E}g(\si
Z)$ is non-decreasing in $\si >0$ for $Z\sim N(0,1)$ and any convex
function $g$. One way to verify the latter statement is as follows. It
suffices to consider the functions of the form $g(u)=(u-t)_{+}$ for real $t$%
; cf. identity () in Pinelis (1994). But the derivative of $\mathsf{E}%
(\si Z-t)_{+}$ in $\si >0$ is $\varphi (t/\si )>0$. Alternatively,
one can prove that $\mathsf{E}g(\si Z)$ is non-decreasing in $\si >0$
by an application of Jensen's inequality. $\Bigr)$ Moreover, by rescaling,
w.l.o.g. $d-c=2$. In other words, then one has the following: 
\begin{equation*}
X=%
\begin{cases}
2r & \text{with probability }1-r, \\ 
2r-2 & \text{with probability }r,%
\end{cases}%
\end{equation*}%
for some $r\in \lbrack 0,1]$. At that, 
\begin{equation*}
Y\sim N(0,1).
\end{equation*}%
Now the right-hand side of inequality \eqref{eq:bern-moment-compar} can be
written as 
\begin{equation}
\E f_{t}(Y)=R(t):=P(t)\varphi (t)-Q(t)\Psi (t),  \label{eq:R(t)}
\end{equation}%
where 
\begin{equation*}
P(t):=8+9t^{2}+t^{4}\quad \text{and}\quad Q(t):=t(15+10t^{2}+t^{4}),
\end{equation*}%
and its left-hand side as 
\begin{equation}
{\E f_{t}(X)=L(r,t):=r(2r-2-t)_{+}^{5}+(1-r)(2r-t)_{+}^{5},}
\label{eq:L(r,t)}
\end{equation}%
so that \eqref{eq:bern-moment-compar} is reduced to the inequality 
\begin{equation}
{L(r,t)\leq R(t)}  \label{eq:L<R}
\end{equation}%
for all $r\in \lbrack 0,1]$ and all real $t$.

Note that \eqref{eq:L<R} is trivial for $t\ge2r$, because then $L(r,t)=0$.

Therefore, it remains to consider two cases: $(r,t)\in B$ and $(r,t)\in C$,
where 
\begin{multline*}
B:=\{(r,t)\colon 0\le r\le1,t\le2r-2\}\quad\text{and}\quad \\
C:=\{(r,t)\colon 0\le r\le1,2r-2\le t\le2r\}.
\end{multline*}

\textbf{Case 1} $(r,t)\in B$.\quad Note that in this case $t\le0$ and, by %
\eqref{eq:L(r,t)}, 
\begin{equation*}
L(r,t)=r(2r-2-t)^5 + (1-r)(2r-t)^5.
\end{equation*}
For $t\ne0$, one has the identity 
\begin{equation}  \label{eq:identity}
{\ \frac{Q(t)^2}{\varphi(t)}\;\partial_t \left(\frac{R(t)-L(r,t)}{Q(t)}%
\right) =Q_2(r,t):=\frac{Q_1(r,t)}{\varphi(t)}-120, }
\end{equation}
where 
\begin{equation*}
Q_1(r,t):=Q^{\prime}(t)L(r,t)-Q(t)\,\partial_t L(r,t),
\end{equation*}
which is a polynomial in $r$ and $t$. Note that 

\begin{equation*}
\partial_r Q_2(r,t) =\frac{\partial_r Q_1(r,t)}{\varphi(t)} \quad\text{and}%
\quad \partial_t Q_2(r,t)=\frac{20\,Q(t)}{\varphi(t)}d(r,t), 
\end{equation*}
where 
\begin{equation*}
d(r,t):=\frac{tQ_1(t)+\partial_t Q_1(t)}{20\,Q(t)}
\end{equation*}
is a polynomial in $r$ and $t$, of degree 2 in $r$. Therefore, the critical
points of $Q_2$ in the interior $\mathrm{int}\,{B}$ of domain $B$ are the
solutions $(r,t)$ of the system of polynomial equations 
\begin{equation*}
\begin{cases}
d(r,t) & =0, \\ 
\partial_r Q_1(r,t) & =0. 
\end{cases}%
\end{equation*}
Further, one has 
\begin{equation*}
\text{$d(r,t)=0$ if and only if $r=r_1(u)$ or $r=r_2(u)$,}
\end{equation*}
where 
\begin{equation*}
u:=2-r-t>0,\quad r_1(u):=\frac{1+u/2}{1+u}\in(0,1),\quad\text{and}\quad
r_2(u):=\frac{2+2u+u^2/2}{2+2u+u^2}\in(0,1).
\end{equation*}

Using the Sturm theorem or the convenient command \textbf{Reduce} of
Mathematica 5.0, one can see that the only solution $u=u_1>0$ of the
algebraic equation $\partial_r Q_1(r,t)|_{r=r_1(u),t=2-r_1(u)-u}=0$ is $%
0.269\dots$, and \\
$Q_2(r,t)|_{r=r_1(u_1),t=2-r_1(u_1)-u_1}<0$. 
As for the
equation \\
$\partial_r Q_1(r,t)|_{r=r_2(u),t=2-r_2(u)-u}=0$, it has no
solutions $u>0$.

Thus, $Q_2<0$ at the only critical point $(r,t)=\bigl(r_1(u_1),2-r_1(u_1)-u_1%
\bigr)$ of $Q_2$ in $\mathrm{int}\,{B}$.

Next, with $u>0$, 
\begin{equation*}
Q_2(r,t)|_{r=0,t=2r-2-u}= -20\left(6+\frac{(2 + u)^5}{\varphi(2+u)} \, 
\left( 7 + 4\,u + u^2 \right)\right)<0.
\end{equation*}
Similarly, with $u>0$, 
\begin{equation*}
Q_2(r,t)|_{r=1,t=2r-2-u}= -20\,\left( 6 + \frac{u^5(3+u^2)}{\varphi(u)}
\right)<0.
\end{equation*}

Now consider the function 
\begin{equation*}
q_2(r):=Q_2(r,t)|_{t=2r-2}.
\end{equation*}

Then $\varphi(2r-2)q^{\prime}_2(r)$ is a polynomial, whose only root $%
r=r_3\in(0,1)$ is $0.865\dots$ . But $q_2(r_3)<0$. Therefore, $Q_2<0$ at the
only critical point of $Q_2$ in the relative interior of the boundary $t=2r-2
$ of domain $B$.

Thus, as far as the sign of $Q_2$ on $B$ is concerned, it remains to
consider the behavior of $Q_2$ as $t\to-\infty$, which is as follows: $%
Q_2(r,t)\sim20(2r-1)^2 t^7\to-\infty<0$ for every $r\ne1/2$ and $%
Q_2(r,t)\sim40t^3(5+t^2)\to-\infty<0$ for $r=1/2$.

(As usual, $a\sim b$ means $a/b\to1$.)

We conclude that $Q_2<0$ on $B$. Hence, in view of \eqref{eq:identity}, the
ratio $\frac{R(t)-L(r,t)}{Q(t)}$ is decreasing in $t$ on $B$.

Next, note that $\varphi(t)$ and $1-\Psi(t)$ are $o(1/|t|^p)$ for every $p>0$
as $t\to-\infty$. Hence, in view of \eqref{eq:R(t)}, one has the following
as $t\to-\infty$: $R(t)-L(r,t)=-Q(t)-L(r,t)+o(1)\sim-10(2r-1)^2t^3\to\infty$
for every $r\ne1/2$ and $R(t)-L(r,t)=-10t\to\infty$ for $r=1/2$.

Hence, $\frac{R(t)-L(r,t)}{Q(t)}<0$ for each $r\in(0,1)$ and all $t<0$ with
large enough $|t|$. Since $\frac{R(t)-L(r,t)}{Q(t)}$ is decreasing in $t$ on 
$B$, one has $\frac{R(t)-L(r,t)}{Q(t)}<0$ on $B$, whence $L(r,t)\le R(t)$ on 
$B$ (because $Q(t)\le0$ on $B$).

It remains to consider

\textbf{Case 2} $(r,t)\in C$.\quad Here, letting $v:=2r-t$, one has $0\le v\le2$, and, by \eqref{eq:L(r,t)}, 
\begin{equation*}
L(r,t)=(1-r)(2r-t)^5.
\end{equation*}
Let us use here notation introduced in the above consideration of Case 1.
Then 
\begin{equation*}
d(r,t)|_{t=2r-v} = -(1-r)v^3\left(1-\frac r2 v\right)<0
\end{equation*}
for $(r,t)=(r,2r-v)\in\mathrm{int}\, C$. This implies that $Q_2$ has no
critical points in $\mathrm{int}\, C$.

Next, with $v>0$, 
\begin{equation*}
Q_2(r,t)|_{r=0,t=2r-v}=  -20\left(6 +\frac{v^5(3+v^2)}{\varphi(t)}\right)<0.
\end{equation*}

On the boundaries $r=1$ and $t=2r$ of $C$, one has $Q_2=-120<0$. The
boundary $t=2r-2$ of $C$ is common with $B$, and it was shown above that $%
Q_2<0$ on that boundary as well.

Thus, $Q_2<0$ on $C$. Since $Q(t)=0$ only for $t=0$, it follows that the
ratio $\frac{R(t)-L(r,t)}{Q(t)}$ is decreasing in $t$ on $C$.

Hence, just as on $B$, one has that $L(r,t)<R(t)$ on $C_-:=\{(r,t)\in
C\colon t\le0\}$.

Moreover, $\frac{R(t)-L(r,t)}{Q(t)}=\frac{R(t)}{Q(t)}>0$ for $t=2r$, since $%
Q>0$ on $C_+:=C\setminus C_-=\{(r,t)\in C\colon t>0\}$. Because $\frac{%
R(t)-L(r,t)}{Q(t)}$ is decreasing in $t$, one has $\frac{R(t)-L(r,t)}{Q(t)}>0
$ on $C_+$ and hence $L(r,t)<R(t)$ on $C_+$.
\end{proof}

\begin{proof}[Proof of Theorem \ref{th:bernoulli-improved}]
This proof is similar to the proof of Theorem~2.1 in \cite{binom} and Theorem~\ref{th:bernoulli}, but based on the following lemma, instead of Lemma~3.2 in \cite{binom} or \ref{lem:bernoulli}. (As in the proof of Theorem \ref{th:bernoulli}, here one has also to refer to Lemma~3.1 in \cite{binom}, rather than Lemma \ref{lem:reduction1}.)
\end{proof}

\begin{lemma}
\label{lem:bernoulli-si} Suppose that $X$ is a r.v. such that $\mathsf{E} X=0$, $X\le d$ with probability 1, and 
$\E X^2\le\si^2$, for some positive constants $d$ and $\si$. Let
$$s:=\frac12\left( d+\frac{\si^2}d \right).$$
Let $Z\sim N(0,1)$.
Then for all $f\in\G^{(5)}$ 
\begin{equation}  \label{eq:bern-moment-compar-si}
\E f(X)\le\E f(sZ).
\end{equation}
\end{lemma}

\begin{proof}
In view of \eqref{eq:F-al-beta}, one has $\F^{(5)}\subseteq\F^{(2)}$. Therefore, by Lemma~3.2 in \cite{binom},
one may assume without loss of generality that here $X=d\cdot X_a$, where $a=\si^2/d^2$. Now it is seen that Lemma \ref{lem:bernoulli-si} follows from Lemma \ref{lem:bernoulli}.
\end{proof}

\subsection{Proofs for Section \ref{max}} \label{subsec:max}

\begin{proof}[Proof of Theorem \ref{th:max}]
Lemma~\ref{lem:reduction1} and Lemma~3.1 in \cite{binom}
reduce Theorem \ref{th:max} to the case when $(S_n)$ is a martingale, and then Theorem \ref{th:max} follows by Doob's inequality \eqref{eq:doob-0}. 
\end{proof}

\begin{proof}[Proof of Theorem \ref{th:al-beta-Mn}]
For every $y>t$, by Doob's inequality,
$$\P(M_n\ge y)\le \frac{\E(S_n-t)_+ I\{M_n\ge y\} }{y-t} .$$
Hence, letting
\begin{equation}\label{eq:J}
J(u):=\int_x^u \frac{\beta(y-x)^{\beta-1}}{y-t}\,\d y\ I\{u>x\} \quad\text{and}\quad
\al':=\frac\al{\al-1},
\end{equation}
and using Fubini's theorem, one has
\begin{align}
\E(M_n-x)_+^\beta &= \int_x^\infty \beta(y-x)^{\beta-1} \P(M_n\ge y) \,\d y \notag \\
&\le \int_x^\infty \beta(y-x)^{\beta-1} \frac{\E(S_n-t)_+ I\{M_n\ge y\} }{y-t} \,\d y \notag \\
&=\E \int_x^\infty \beta(y-x)^{\beta-1} \frac{(S_n-t)_+ I\{M_n\ge y\} }{y-t} \,\d y \notag \\
&= \E(S_n-t)_+ J(M_n) \notag \\
&\le \left( \E(S_n-t)_+^\al \right)^{1/\al}\, \left( \E J(M_n)^{\al'} \right)^{1/\al'}, \label{eq: le J}
\end{align}
by H\"older's inequality. 

Observe that for all real $u$
\begin{equation}\label{eq:J-le}
J(u) \le c^{1/\al} (u-x)_+^{\beta/\al'}, \quad\text{where}\quad c:=\frac{k_{1;\al,\beta}}{(x-t)^{\al-\beta}}.
\end{equation}
Indeed, introducing new variables $\si:=\frac{u-x}{x-t}$ and $s:=\frac{y-x}{x-t}$, one can see that, for $u>x$,
\begin{align*}
J(u) &= (x-t)^{\beta-1} \int_0^\si \frac{\beta s^{\beta-1}\,\d s}{1+s} \quad\text{and}\quad \\
c^{1/\al} (u-x)_+^{\beta/\al'} &= k_{1;\al,\beta}^{1/\al} \si^{\beta(1-1/\al)} (x-t)^{\beta-1},
\end{align*}
so that \eqref{eq:J-le} follows, in view of \eqref{eq:k1}.

Now \eqref{eq: le J} and \eqref{eq:J-le} imply \eqref{eq:al-beta-Sn}.
\end{proof}

\begin{proof}[Proof of Theorem \ref{th:max-beta}]
This is similar to the proof Theorem \ref{th:max}, but relies on inequality \eqref{eq:al-beta-Mn} in place of Doob's inequality \eqref{eq:doob-0}. 
\end{proof}

\begin{proof}[Proof of Proposition \ref{prop:k1-k2}]
Introduce
\begin{align*}
f(\si,\al,\beta,\g) &:=
\si^{-\beta(\al-\g)/\g}
\left(\int_0^\si \frac{\beta s^{\beta-1}\,ds}{(1+s)^\g}\right)^{\al/\g}, \\
K(\al,\beta,\g) &:= \sup_{\si>0} f(\si,\al,\beta,\g).
\end{align*}
Then $\si^{-\beta/\al}f(\si,\al,\beta,\g)^{1/\al}=(\E Y^\g)^{1/\g}$, where $Y:=\frac1{1+S}$ and $S$ is a r.v. with density 
$s\mapsto\si^{-\beta}\beta s^{\beta-1}I\{0<s<\si\}.$
Hence, $f(\si,\al,\beta,\g)$ is non-decreasing in $\g$, and then so is $K(\al,\beta,\d)$. 
Therefore,
$$k_{1;\al,\beta}=K(\al,\beta,1)\le K(\al,\beta,\al)=k_{2;\al,\beta}.$$
\end{proof}

\begin{proof}[Proof of Proposition \ref{prop:k1-k3}]
By \eqref{eq:k1},
\begin{equation}\label{eq:k1-r}
k_{1;\al,\beta}=\sup_{\si>0}r(\si)^\al,
\end{equation}
where
$$r(\si):=\frac{f(\si)}{g(\si)},\quad f(\si):=\int_0^\si\frac{\beta s^{\beta-1}\,ds}{1+s},
\quad\text{and}\quad g(s):=\si^{\beta(1-1/\al)}.$$
Note that the monotonicity pattern of 
\begin{equation}\label{eq:r1}
r_1(\si):=\frac{f'(\si)}{g'(\si)}=\frac\al{\al-1}\frac{\si^{\beta/\al}}{1+\si}
\end{equation}
on $(0,\infty)$ is $\uu\se$; that is, there exists some $\si_1(\al,\beta)\in(0,\infty)$ such that $r_1\uu$ (is increasing) on $(0,\si_1(\al,\beta))$ and $r_1\se$ (is decreasing) on $(\si_1(\al,\beta),\infty)$; namely, here 
\begin{equation}\label{eq:si1}
\si_1(\al,\beta)=\frac\beta{\al-\beta}.
\end{equation}
Also, $gg'>0$ on $(0,\infty)$. Hence, it follows from \cite[Proposition 1.9]{waves} that $r$ has one of these monotonicity patterns on $(0,\infty)$: $\uu$ or $\se$ or $\uu\se$ or $\se\uu$ or $\se\uu\se$. However, $r(\si)$ is positive on $(0,\infty)$ and converges to 0 when $\si\downarrow0$ as well as when $\si\to\infty$. This leaves only one possible pattern for $r$: $\uu\se$.   
Hence, there is some $\si(\al,\beta)\in(0,\infty)$, at which $r$ attains its maximum on $(0,\infty)$; moreover, $r'(\si(\al,\beta))=0$, which is equivalent to $r(\si(\al,\beta))=r_1(\si(\al,\beta))$. Thus,
\begin{multline*}
k_{1;\al,\beta}=\sup_{\si>0}r(\si)^\al
=r(\si(\al,\beta))^\al=r_1(\si(\al,\beta))^\al
\le\sup_{\si>0}r_1(\si)^\al \\
=r_1(\si_1(\al,\beta))^\al=k_{3;\al,\beta},
\end{multline*}
in view of \eqref{eq:r1}, \eqref{eq:si1}, and \eqref{eq:k1-k3}.
\end{proof}

\begin{proof}[Proof of Proposition \ref{prop:k1=k3}] 
In the case $\beta=\al>1$, the function $r_1$ given by 
\eqref{eq:r1} is increasing on $(0,\infty)$ to 
$r_1(\infty)=\frac\al{\al-1}$. Hence, so does $r$, 
according to the mentioned \cite[Proposition 1.9]{waves}. 
Now 
Proposition \ref{prop:k1=k3} follows in view of \eqref{eq:k1-r}. 
\end{proof}

\begin{proof}[Proof of Proposition \ref{prop:al-beta}]
Elementary calculus; the optimal value of $u$, when inequality \eqref{eq:al-beta} turns into an equality, is
\begin{equation}\label{eq:u*}
u_*:=\frac{\al x-\beta t}{\al-\beta}>x.
\end{equation}
\end{proof}

\begin{proof}[Proof of Proposition \ref{prop:al-beta-Sn}]
Only that $k_{\al,\beta}$ is the best constant factor needs to be proved. Without loss of generality, $x>0$. Suppose that \eqref{eq:al-beta-Sn} holds with some constant $\tilde k$ in place of $k_{\al,\beta}$; then, by continuity, it holds for the continuous-time martingale $S_v:=B_{v\wedge\tau}$ in place of $S_n$, where $B(\cdot)$ is a standard Brownian motion, $v\ge0$, and 
$$\tau:=\inf\{v\ge0\colon B_v=u_*\text{ or }B_v=\tilde t\};$$
here, $u_*$ is defined by \eqref{eq:u*} and $\tilde t:=(-1)\wedge t$. Note that $\E\tau=u_*|\tilde t|$ and $p:=\frac{|\tilde t|}{|\tilde t|+u_*}>0$. It follows that
$$p\cdot(u_*-x)^\beta=
\E(S_\infty-x)_+^\beta\le\tilde k\frac{\E(S_\infty-t)_+^\al}{(x-t)^{\al-\beta}}
=\tilde k p\frac{(u_*-t)_+^\al}{(x-t)^{\al-\beta}}.
$$
Because $k_{\al,\beta}$ is the best constant in \eqref{eq:al-beta}, it follows now that $\tilde k\ge k_{\al,\beta}$.
\end{proof}

\subsection{Proofs for Section \ref{concentr}} \label{subsec:concentr}

The proofs here are based on the improvements given in \cite{pin81a} and \cite{pin-sakh} of the method of Yurinskii (1974) \cite{yurin}; cf. \cite{mcdiarmid1,mcdiarmid2} and \cite{bent-isr}. 

For a r.v. $Y$ as in Theorem \ref{th:concentr}, consider the martingale expansion 
$$ Y - \E Y = \xi_1+\dots+\xi_n,$$
of $Y - \E Y$
with the martingale-differences 
\begin{equation}\label{eq:xi_i}
\xi_i:= \E_i Y - \E_{i-1} Y.
\end{equation}
where $\E_i$ denotes the conditional expectation given $H_{\le i}:=(X_1,\dots,X_i)$.
For each $i$ pick an arbitrary non-random $x_i$, and introduce the r.v.
\begin{equation}\label{eq:eta_i}
\eta_i := Y - \tilde Y_i ,\quad\text{where}\quad
\tilde Y_i := g(X_1,\dots,X_{i-1},x_i,X_{i+1},\dots,X_n).
\end{equation}

\begin{proof}[Proof of Theorem \ref{th:concentr} and Proposition \ref{prop:LipRelaxed}]
Nore that, for the function $\Xi_i$ defined by \eqref{eq:Xi_i}, one has $\Xi_i(X_1,\dots,X_i)=\xi_i$, where $\xi_i$ is defined by \eqref{eq:xi_i}. 
It follows from \eqref{eq:xi_i} that
\begin{equation}\label{eq:CD2}
C_{2,i-1} \le \xi_i \le D_{2,i-1}\quad\text{and}\quad
D_{2,i-1} - C_{2,i-1} \le 2\hat r_i\le2 r_i,
\end{equation}
where $r_i$ and $\hat r_i$ are given by \eqref{eq:r_i} and \eqref{eq:LipRelaxed}, and
\begin{align*}
C_{2,i-1} &:= \inf_{x_i} \E_{i-1} (-\eta_i) = \inf_{x_i} \E_{i-1} \tilde Y_i -\E_{i-1} Y\quad\text{and} \\
D_{2,i-1} &:= \sup_{x_i} \E_{i-1} (-\eta_i) = \sup_{x_i} \E_{i-1} \tilde Y_i -\E_{i-1} Y
\end{align*}
are 
$H_{\le(i-1)}$-measurable.  
Now Proposition \ref{prop:LipRelaxed} -- and hence Theorem \ref{th:concentr} -- follow by Theorem \ref{th:bernoulli} 
and Corollary \ref{cor:bern-beta}.
\end{proof}

\begin{proof}[Proof of Theorem \ref{th:concentr-modif}]
This proof is similar to that of Theorem \ref{th:concentr} and Proposition \ref{prop:LipRelaxed}, but based on Theorem \ref{th:bernoulli-improved} in place of Theorem \ref{th:bernoulli} 
and Corollary \ref{cor:bern-beta}.
(Note that $\E\Xi_i(x_1,\dots,x_{i-1},X_i)^2$ is the same as conditional expectation $\E_{i-1}\xi_i^2$ given that $X_1=x_1,\dots,X_{i-1}=x_{i-1}$.) 
\end{proof}

\begin{proof}[Proof of Proposition \ref{prop:sigma-xi-bound}]
For each $i$,
\begin{equation}\label{eq:xi-eta}
\xi_i = \E_i \eta_i - \E_{i-1} \eta_i, 
\end{equation}
because $ \E_i \tilde Y_i = \E_{i-1} \tilde Y_i$, in view of the independence of the $X_i$'s. 
Hence and by \eqref{eq:Lip}, for any given $x_i$,
\begin{equation}\label{eq:|eta|}
|\eta_i| \le \rho_i(X_i,x_i) 
\end{equation}
with probability 1.
It follows from \eqref{eq:xi-eta} and \eqref{eq:|eta|} that, for any $x_i$,
\begin{multline*}
\E_{i-1}\xi_i^2 = \E_{i-1}(\E_i \eta_i - \E_{i-1} \eta_i)^2
=\Var_{i-1}(\E_i\eta_i) \le \E_{i-1}(\E_i \eta_i)^2 \le \E_{i-1}\E_i\eta_i^2 \\ 
= \E_{i-1}\eta_i^2 
\le \E_{i-1}\rho_i(X_i,x_i)^2 = \E\rho_i(X_i,x_i)^2, 
\end{multline*}
which proves \eqref{eq:sigma-xi-bound}; here, $\Var_{i-1}$ denotes the conditional variance given $H_{i-1}$. 

To prove \eqref{eq:sigma-bound}, suppose in addition that the function $g$ is convex in each of its arguments, as stated in the second part of Proposition \ref{prop:sigma-xi-bound}. Let $\tilde\E_i$ denote the conditional expectation given $X_1,\dots,X_{i-1},X_{i+1},\dots,X_n$. Then, for all $i$, by Jensen's inequality,
\begin{multline*}
\E_{i-1}Y=\E_{i-1}\tilde\E_i Y =\E_{i-1}\tilde\E_i g(X_1,\dots,X_n) \\
\ge \E_{i-1}g(X_1,\dots,X_{i-1},\tilde\E_i X_i,X_{i+1},\dots,X_n) \\
= \E_{i-1}g(X_1,\dots,X_{i-1},\E X_i,X_{i+1},\dots,X_n)
= \E_{i-1}\tilde Y,
\end{multline*}
in view of \eqref{eq:eta_i},
if $x_i$ is chosen to coincide with $\E X_i$; hence, 
$$\E_{i-1}\eta_i=\E_{i-1}Y-\E_{i-1}\tilde Y\ge0.$$
This and formulas \eqref{eq:xi-eta} and \eqref{eq:|eta|} imply that
$$ \xi_i \le \E_i \eta_i \le \rho_i(X_i,\E X_i),$$
which is equivalent to \eqref{eq:sigma-bound}.
\end{proof}

\begin{proof}[Proof of Corollary \ref{cor:banach}]
This follows immediately from Theorem \ref{th:concentr-modif} and Proposition \ref{prop:sigma-xi-bound}, with $\rho_i(\tilde x_i,x_i)=\|\tilde x_i-x_i\|$.
\end{proof}

\end{document}